\documentclass[12pt]{amsart}
\usepackage{amsfonts}
\usepackage{amsmath, amssymb}
\usepackage{mathtools, amsthm}
\usepackage{graphicx}
\usepackage{tikz}
\usepackage{mathdots}
\usepackage{yhmath}
\usepackage{cancel}
\usepackage{color}
\usepackage{siunitx}
\usepackage{array}
\usepackage{multirow}
\usepackage{bigstrut}
\usepackage{bigdelim}
\usepackage{gensymb}
\usepackage{tabularx}
\usepackage{multirow} 
\usepackage{booktabs}
\usetikzlibrary{fadings}
\usetikzlibrary{patterns}
\newtheorem{theorem}{Theorem}
\newtheorem{lemma}[theorem]{Lemma}

\theoremstyle{definition}

\theoremstyle{remark}
\newtheorem{remark}[theorem]{Remark}

\numberwithin{equation}{section}
\numberwithin{theorem}{section}

\newtheoremstyle{case}{}{}{}{}{}{:}{ }{}
\theoremstyle{case}

\newcommand{\zpstar}{\ensuremath{\mathbb{Z}_p^*}}
\newcommand{\zp}{\ensuremath{\mathbb{Z}_p}}
\begin{document}

\title[Solutions of the equation $x^{n} + y^{n} = z^{n}$ in the finite fields $\zp$]{On The solutions of the Diophantine equation  $x^n + y^n = z^n$ in the finite fields $\zp$.}


\author{Silvia R. Valdes}
\address{Department of Mathematics and Statistics \\
University of Southern Maine \\
Portland, ME 04104-9300}
\curraddr{ }
\email{svaldes@maine.edu}
\thanks{}

\author{Yelena Shvets}
\address{Department of Mathematics and Statistics \\
University of Southern Maine \\
Portland, ME 04104 -9300 }
\curraddr{}
\email{yelena.shvets@maine.edu}
\thanks{}

\subjclass[2010]{Primary }


\date{\today}

\dedicatory{}

\begin{abstract}
Let $p$ be a prime integer, $\zp$  the finite field  of order $p$ and $\zpstar$ its multiplicative cyclic group. We consider the Diophantine equation $x^n + y^n = z^n$ with $1 \leq n \leq \frac{p - 1}{2}$. \\
Our main aim in this paper is to give best-possible conditions or relationships between the exponent $n$ and the prime $p$ to determine the existence of nontrivial solutions of the diophantine equation $x^n + y^n = z^n$ with $1 \leq n \leq p -1 $, in finite fields $\zp$. 
\end{abstract}

\maketitle

\section{Problem Statement} 
Let $p$ be a prime and let $\zp$ be the finite field with $p$ elements. Let's denote by $\zpstar$ the multiplicative cyclic group of order $p-1$, comprised of the non-zero elements in $\zp$.\\

The purpose of this work is to give a means of algorithmically constructing all solutions to
\begin{equation} \label{eq:BasicM}
x^n + y^n \equiv z^n\pmod p,
\end{equation}

or equivalently

\begin{equation} \label{eq:BasicZ}
x^n + y^n = z^n, \text{ in $\zpstar$}
\end{equation}
for all integer exponents $n$.  The cyclic nature of the group $\zpstar$ allows us to restrict our consideration first to $1 \leq n \leq p -1$ and eventually to $1 \leq n \leq \frac{p -1}{2}$  .\\

More precisely, we develop a procedure which allows us, for any given prime $p$, to  enumerate the powers $1 \leq n \leq p -1$, for which the equation has non-trivial solutions, and for which it does not. It also enables us to construct all such solutions, provided that a generator of $\zpstar$ has been identified.\\

We start with an observation that given any solution $(x_0,y_0,z_0)$ of ($\ref{eq:BasicM}$), and any unit $u \in \zpstar$, a triple $(u\cdot x_0,u\cdot y_0,u \cdot z_0)$ is also a solution,  to which we refer as a $u$-multiple.\\
\\
With this in mind, we observe that the following triples are always solutions, for any positive integer $n$ and any prime $p$:
\begin{itemize}
\item (0,0,0),
\item (1,0,1),
\item (0,1,1).
\end{itemize}
We call such solutions and their $u$-multiples {\em trivial solutions}.\\

If there is a solution of the form $(x_0, y_0, 0)$, with $x_0, y_0  \in \zpstar$, then we call this solution and its $u$-multiples,  {\em type-0 solutions}.\\
 
Similarly, given any solution of the form $(x_0, y_0, z_0)$, with $x_0, y_0 , z_0 \in \zpstar$, we call this solution and its $u$-multiples, {\em type-1 solutions}. The name is motivated by the observations that all such solutions are $u$-multiples of a solution of the form  $(z_0^{-1}x_0, z_0^{-1}y_0, 1)$.\\
 
 \section{General results}\label{sec:GR}
 Here we provide a few general remarks that are relevant for both types of non-trivial solutions.
 Observe that if $p=2$ and $n = p-1$, then $(1,1,0)$ is a type-0  solution of (\ref{eq:BasicM}).
 \begin{remark} [n=1 and n=p] \label{rem:1}
 For any prime $p$, the equation (\ref{eq:BasicM}) has type-0 solutions for $n = 1$ and $n=p$. This is an immediate consequence of Fermat's Little Theorem:
 $$a^p \equiv a\pmod p,$$
 where $a$ is any element in  $\zpstar$.\\
 Furthermore, for any prime $p >2$, the equation (\ref{eq:BasicM}) has type-1 solutions for $n = 1$ and $n=p$. \\
 
 Let $x_0  \neq 1 \in \zpstar$, then 
 $$ 
 x_0 + (p+1 - x_0) \equiv 1\pmod p,
 $$
 and $(x_0, p+1-x_0,1)$ is a type-1 solution of $x+y \equiv z \pmod p$ and $x^p + y^p \equiv  z^p \pmod p$.
 \end{remark}
 \begin{remark} [n=p-1] \label{rem:2}
 For $p>2$ and $n = p-1$, the equation (\ref{eq:BasicM}) has no non-trivial solutions in $\zp$. Which, once again follows from Fermat's Little Theorem:
 $$a^{p-1} \equiv 1\pmod p,$$
 where $a \in \zpstar$.
 \end{remark}
 \begin{lemma}
Let $1 < t < p - 1$ with  $p > 2 $  and $\gcd(t, p - 1) = 1$. Then the set of $t$'th powers in $\zp$ is exactly as the set of the first powers and hence any solution of  $x + y = z$ will give rise to a solution of  $x^t + y^t = z^t$  . \\
\end{lemma}
Proof: \\
\\
Let $0 \neq a$; $0 \neq b$ in $\zp$ with $a \neq b$. It is enough to prove that $a^t \neq b^t$ in $\zp$. Thus the set $\{0, 1, 2^t, 3^t, \cdots ,(p - 1)^t \}$ is exactly the same as the set $\{0, 1, 2, 3, \cdots , (p - 1)\}$.
Assume that $a^t \equiv b^t \pmod p$ , then $(ab^{-1})^t \equiv 1 \pmod p$. So if $\lambda$ is the order of $ab^{-1}$ , $\lambda \mid t$ and $\lambda  \mid (p-1)$; that is, $\lambda \mid \gcd(t,p - 1)$, ie $\lambda = 1$, thus $a = b$, a contradiction.
\\
 \begin{lemma}
 For $k \mid(p - 1)$ the set of the $k$-powers of elements in $\zpstar$, where $p > 2$, constitute the unique cyclic subgroup of $\zpstar$ of order $\frac{p - 1}{k}$.  
 \end{lemma}
 Proof: \\
 \\
 Let $0 \neq a$ in $\zp$ with $|a| = p - 1$; that is, $\zpstar = < a >.$  Thus $a^k$ has order $\frac{p - 1}{k}$. Let $b = a^k$, then 
 $< b > = \{b, b^2, \cdots, b^{\frac{p - 1 - k}{k}}, b^{\frac{p -1}{k}} \}$ the unique subgroup of $\zpstar$ of order $\frac{p - 1}{k}$. Now if $0 \neq x_0$ in $\zp$ with $x_0 \equiv t^k \pmod p$ with $t \in \zpstar$, then $x_0^{\frac{p - 1}{k}} \equiv 1 \pmod p$. Thus $\lambda$,  the order of $x_0$ is a divisor of $\frac{p - 1}{k}$, we write $\frac{p - 1}{k} = \lambda \cdot \delta$ with $\delta \geq 1$. Because $< x_{0} >$ is a subgroup of $< b >$, we have to consider $x_{0} = b^{\delta} = a^{k\delta}$, and the result follows.  \\
 
Recall that in general, if $a$ has order $k$ modulo $p$ and $h > 0$, then $a^h$ has order $\frac{k}{\gcd(h,k)} \pmod p$ ([1] Theorem 8.3). Observe that since $a^{-k} \in < b >$, and $a^{-k} = a^{p - 1 - k}$. The set of the $k$-powers of elements in $\zpstar$ equal the set of all $p-1-k$-powers of elements in $\zpstar$.\\

 With this in mind we see that it is enough to consider the $k$-powers of elements in $\zpstar$ with $1 \leq k < \frac{p - 1}{2}$. The special case when $k = \frac{p - 1}{2}$ will be addressed in section \ref{sec:T1}. \\
 \begin{remark}
 Let q be a prime such that $q^s \mid (p-1)$ and $q^{s + l} \nmid (p - 1)$, for all $l \geq 1$.  In that case the set of the $q^{s + l}$- powers of elements in $\zpstar$ is exactly as the set of the $q^s$-powers of elements in $\zpstar$. \\
 
Let $\zpstar = <a>$. The set of the $q^{s}$- powers of elements in $\zpstar$ is the subgroup $< b >$ with $b = a^{q^{s}}$ of order $\lambda = \frac{p - 1}{q^{s}}$ of $\zpstar$. Now let $c = a^{q^{s + l}}$. \\
$|c| = \frac{p - 1}{\gcd(q^{s + l}, p - 1)} = \frac{p - 1}{\gcd(q^{s + l}, q^{s} \lambda)} = \frac{p - 1}{q^{s}}  = \lambda$.  Since $c \in <b>$  and $|<c>| = \lambda = |<b>| $ we must have $< c > = < b >$. 
 
Similarly, under the same assumptions, it can be seen that if $\gcd(m,p-1) = 1$, the the group of the $q^{s} \times m $-powers of elements in $\zpstar$ is exactly the same as the group of $q^s$ powers.\\
 \end{remark}
 
The following remark provides a generalization.
 
 \begin{remark}
 Let $t \nmid (p - 1)$ with $\gcd(t, p - 1) = d > 1$. We write $t = dl$  with $ 1 < l < p - 1$ and $p - 1 = d s = \frac{t}{l} s$.  \\
 Here we have a couple of cases: 
 \renewcommand{\theenumi}{\roman{enumi}}%
 \begin{enumerate}
 \item If $\gcd(l, p - 1) = 1$, then the set of the $t$-powers of elements of $\zpstar$ is exactly as the set of the $d$-powers of elements of $\zpstar$, which is a subgroup of order $\frac{p - 1}{d}$.
 \item If $\gcd(l, p - 1) = x > 1$, then the set of $t$-powers of elements of $\zpstar$ would be a subset of the subgroup of the $d$-powers and a subset of the subgroup of the $x$-powers.
 \end{enumerate}
 
 For example: Let $p = 23$, and $t = 8$. In this case we are in case (ii) of the previous remark with $d = 2$, $l = 4$. \\

 \end{remark}
\begin{remark}
Let $\zpstar = < a >$ be the cyclic group of the nonzero elements of the finite field $\zp$, and $U(\mathbb{Z}_{p - 1})$ be the group of units of the finite ring $\mathbb{Z}_{p - 1}$. \\
Consider the following function 
\[\zpstar \longrightarrow U(\mathbb{Z}_{p - 1}), \;\; \mbox{defined by} \]
\[ \hspace{-.65in}a \longrightarrow \log_a(a) = 1 \]
\[ b = a^{\mu} \longrightarrow \log_a(b) = \log_a(a^{\mu}) = \mu\]
\noindent
with $b$ any element of $\zpstar$ and $\mu$ the least such power.\\
We also have an inverse mapping:
\[ U(\mathbb{Z}_{p - 1}) \longrightarrow \zpstar \]
\[ 1 \longrightarrow a^1 = a \]
\[\lambda \longrightarrow a^{\lambda} \]
\noindent
\end{remark}

\begin{remark}
Let $p$ be an odd prime, then there are the following cases: \\
\renewcommand{\theenumi}{\roman{enumi}}%
\begin{enumerate}
\item If $p \equiv 1 \pmod 4$, if $a$ is a generator of $\zpstar$, then $-a$ is also a generator. \\
\item  If $p \equiv 3 \pmod 4$, if $a$ is a generator of $\zpstar$, then $-a$ is an element of order $\frac{p - 1}{2}$. That is,  $-a$ generates the subgroup of the 2nd-powers of elements in $\zpstar$. \\
\end{enumerate}
\end{remark}
To address statements (1) and (2) of remark 3.8 (above), consider the following:

Since $a$ is a generator of $\zpstar$, $a^{p - 1} \equiv 1 \pmod{p}$, then $a^{(p - 1)/2} \equiv {-1\pmod{p}}$; that is,  $-a \equiv a^{(p + 1)/2} \pmod{p}$.

In part (1) $ p = 4k + 1$, then $-a \equiv a^{2k + 1} \pmod p$, hence $-a$ has order $\frac{p - 1}{\gcd(2k + 1, p - 1)} = \frac{p - 1}{gcd(2k + 1, 4k)} = p - 1$ ([1], theorem 8.3).

In part (2) $p = 4k + 3$, then $-a \equiv a^{2k + 2} \pmod p$, hence $-a$ has order $\frac{p - 1}{\gcd(2k + 2, 4k + 2)} = \frac{p - 1}{2}$ ([1], theorem 8.3).

\section{Type-0 solutions}\label{sec:T0}
Observe that given a type-0 solution $(x_0,y_0, 0)$, of (\ref{eq:BasicM}), $(x_0, y_0)$ also solves:
\begin{equation}\label{eq:T0M}
x^n + y^n \equiv 0 \pmod p.
\end{equation}\label{eq:T0}

In that case, $y_0^{-1}x_0$ solves:
\begin{equation}
x^n \equiv -1 \pmod p.
\end{equation}
\begin{remark} \label{rem:HalfWay}[$n = \frac{p - 1}{2}$]
For any prime $p$, the equation $x^n + y^n \equiv z^n \pmod p$ has a nontrivial type-0 solution. This is a direct consequence  of Fermat's little Theorem. That is,
\[a^{\frac{p -1}{2}} \equiv  \pm 1 \pmod p\]
\noindent
\end{remark}
\begin{lemma}\label{lemma:T0} 
Let $p$ be an odd prime. Consider the equation $x^{n} + y^{n} \equiv z^n \pmod p$, with $1 \leq  n \leq \frac{p-1}{2}$.  This equation has a type-0 solution if and only if $-1$ is an $n$-th power.
With this in mind,  we have the following cases:\\

\renewcommand{\theenumi}{\roman{enumi}}%
\begin{enumerate}
\item  if $n$ is odd, then there is a type-0 solution;  \\
\item if $n = 2^k$ with $k \geq 1$, then a type-0 solution exists if and only if $\zpstar$ has an element of order $2^{k+1}$;  equivalently, if and only if $p \equiv 1 \pmod{ 2^{k+1}}$;\\
\item if $n = 2^k \cdot t$ with $k \geq 1$ and $t > 1$, an odd integer, then a type-0 solution exists if and only if $\zpstar$ has an element of order $2^{k+1}$;  equivalently, if and only if $p \equiv 1 \pmod {2^{k+1}}$.\\
\end{enumerate}
\end{lemma} 
Proof:\\
By observations above remark (\ref{rem:HalfWay}), we will work with the equation $x^n \equiv -1 \pmod p$. \\
(i) If $n$ is odd, $(p-1)^n = (-1)^n = -1$ then $p-1$ is a desired type-0 solution. Moreover, if $a = w^n$, then the pairs $(w, -w)$ is another such solution.\\
(ii) Let $n = 2^k$,  we know that $x^n \equiv -1 \pmod p$ has a solution if and only if \\
 $(-1)^{(p-1)/gcd(2^k, p-1)} \equiv 1 \pmod{p}$, ([6], page 165). This, in particular, implies that $(p-1)/gcd(2^k, p-1)$ is even. So we can let $p= (2^{k+1} )\cdot \ell + r$, where $r = 1,3,5,\ldots, 2^{k+1}-1$.  \\
 
 We can then write $p -1 = (2^{k+1}) \cdot \ell + (r-1)$. Let $r-1 = 2^{k_1}\cdot t$ for some $1 \leq k_1 \leq k$ and $t$ an odd number or zero. In that case $gcd(2^k, p-1) = 2^{k_1}$.\\ 
 
 Now we write $\frac{p-1}{gcd(2^k, p-1)} = \frac {(2^{k+1}) \cdot \ell + (2^{k_1})\cdot t}{ 2^{k_1}} = (2^{(k-k_1) +1} ) \cdot \ell + t$, which is even. \\
 
 This can happen if and only if $t=0$ and so $r=1$, giving us $$p-1 = 2^{k+1} \cdot \ell \iff p \equiv \pmod{ 2^{k+1}}$$. 
 
 (iii) This follows from (i) and (ii) since for odd $t$,  $$w^{2^k \cdot t} = (w^{2^k} )^t \equiv -1 \pmod {p} \iff (w^{2^k} )\equiv -1 \pmod {p}$$.

Observe that if $n = 2$, we have that following equivalent statements: \\
(i) The equation $x^2 + 1 = 0$ has a solution in $\zpstar$; \\
(ii) $\zpstar$ has an element of order 4; \\
(iii) $p$ is an odd prime with $p \equiv 1 \pmod 4$; \\ 
$\zpstar$ has an element of order 4 if and only if  $4 \mid (p-1)$ if and only $p - 1 = 4t$ for some integer $t$ if and only if  $p \equiv 1 \pmod 4$.

\begin{remark} 
Lemma \ref{lemma:T0} can be restated as: Let $p$ be an odd prime. For $n$ even, the equation $x^n + 1 = 0$,  has a solution in $\zpstar$, if and only if the cyclic group $\zpstar$ has an element of order $2n$. \\ 
\end{remark} 
\begin{lemma} The equation $x^n = -1$ has a solution  in $\zpstar = < a >$ if and only if $0_{\zp} \in \{ < b > + 1\}$ where $< b >$ is the cyclic subgroup of $\zpstar$ of order $\frac{p - 1}{n}$.
\end{lemma} 
Proof:\\
Let $w_0$ be  a solution of the equation $x^n = -1$ in $\zpstar$, ie $w_0 + 1 = 0$; that is $w_0 + 1 \in \{ < b > + 1 \}$ with $w _0 \in < b >$. By the other hand, if $0 \in \{ < b > + 1 \}$, then $0 = b^t + 1$ with $b^t \in < b >$ ($1 \leq t \leq \frac{p - 1}{n}$ and $b = a^n$); then $0 = (a^t)^n + 1$ with $a^t \in \zpstar$. That is, $a^t$ is the desire solution. \\
\\
We combine the above results in the following theorem: \\
\begin{theorem}
let $p$ be an odd prime. Let $\zpstar$ be the multiplicative cyclic group of the finite field $\zp$. Let $n  \mid (p-1)$, with $n$ even, then the following statements are equivalent. \\
(1) The equation $x^n + y^n = 0$ has a nontrivial type-0 solution. \\
(2) The cyclic group $\zpstar$ has an element of order $2n$. \\
(3) $ 0_{\zp} \in \{ < b > + 1 \}$, where $< b >$ is the cyclic subgroup of $\zpstar$ of order $\frac{p - 1}{n}$. \\
\end{theorem}

\section{Type-1 solutions}\label{sec:T1}
Observe that given a type-1 solution $(x_0,y_0,z_0)$ of our initial diophantine equation $x^n + y^n = z^n$,  we have that $y_0^{-1} x_n$ solves 
\begin{equation}
x^n + 1 = z^n. 
\end{equation}
Conversely, any solution $(x_0, z_0)$ of the last equation will produce a solution of original equation, of the form with $(y_0 x_0, y_0, y_0 z_0)$ where $y_0$ is any element of $\zpstar$.
\begin{remark}
Alternatively, the solutions of the initial diophantine equation $x^n + y^n = z^n$ are in one-to-one correspondence with the solutions of: 
\begin{equation}
x^n + y^n = 1
\end{equation}
by using  $z_0^{-1}x_0$.  
\end{remark}
From the results exhibited in section \ref{sec:GR} we know that our work will be focused on the existence of solutions of $x^n + 1 = z^n$ (or $x^n + y^n = 1$) for $n \mid (p-1)$ with $1 < n < p - 1$. \\
We observe the following: \\
(1) Let $1 < n < p - 1$, and let $(u,v)$ be a solution of $1 + x^{n} = z^{n}$. \\
If $n$ is even, then $(-u, -v)$ is a solution of $1 + x^{n} = z^{n}$. \\
If $n$ is odd, then $(p -v, p - u)$ is a solution of $1 + x^{n} = z^{n}$. \\
\\
(2) Let $n = \frac{p - 1}{2} - k$ and $t = \frac{p - 1}{2} + k$, then $(u,v)$ is a solution of $1 + x^{n} = z^{n}$ if and only if $(u^{-1},v^{-1})$ is a solution of $1 + x^{t} = z^{t}$. This observation is also true for the Type-0 solutions.  \\
\\
Therefore, considering the symmetries of the group $\zpstar$, it is enough to focus our study on the existence of the solutions of the equations described above when $1 < n < \frac{p - 1}{2}$. \\
From the results in section \ref{sec:GR}, we already know when the diophantine equation has type-0 solutions. Now we exhibit an algorithm, which is going to give us the exponent(s) for which the diophantine equation has type-1 solutions, and a way to generate such solutions. This algorithm, which we call ``next in line", is described below.  \\
\\
Let $a$ be any generator of $\zpstar$, ie $\zpstar = < a >$. It's worth remarking that among the 78498 odd primes up to $10^6$, the cyclic group $\zpstar$ has a generator less than or equal to 6 ([1], page 156). \\  Raise the generator to each of the powers 1 to $p-1$, modulo $p$, and sort the resulting array:
\[\left. \begin{array}{cccccccc}
  1 & 2 & \cdots & \frac{p - 1}{2} & \frac{p + 1}{2} & \cdots & p - 2 & p - 1 \\
  {} & {} & {} & {} & {} & {} & {} & {} \\
 \Big{\updownarrow} & \Big{\updownarrow} & \cdots & \Big{\updownarrow} & \Big{\updownarrow}  & \cdots  & \Big{\updownarrow} & \Big{\updownarrow} \\ 
  {} & {} & {} & {} & {} & {} & {} & {} \\
 a^{{\alpha}_{1}} & a^{{\alpha}_{2}} & \cdots  & a^{{\alpha}_{\frac{p - 1}{2}}} & a^{{\alpha}_{\frac{p + 1}{2}}} & \cdots & a^{{\alpha}_{p - 2}} & a^{{\alpha}_{p - 1}} \\
 \end{array} \right. \]
 \noindent 
 \\
with $\alpha_{i} \in \{1, 2, 3, \cdots, p - 1 \}$. Also, we
observe that: $\alpha_{1} = p - 1$ and $\alpha_{p - 1} = \frac{p - 1}{2}$. \\
\\
Consider $x = a^{\alpha}$ and $y = x + 1 = a^{\beta}$; that is, $x$ and $y$ are two consecutive elements in $\zp$. Notice also that $a^{\alpha} \equiv a^{ \alpha + k(p-1)} \pmod p,$ for any integer $k$.\\
\\
Let $d = \gcd(\alpha , \beta)$; we write $\alpha= dt_{a}$ and $\beta = dt_{b}$ with $t_{a}$ and $t_{b}$ integers and $\gcd(t_{a}, t_{b}) = 1$. 
\[a^{\beta} = a^{\alpha} + 1,  \text{ in }\zp, \mbox{ can be written as} \]
\[(a^{t_{b}})^{d} = (a^{t_{a}})^{d} + 1,  \text{ in }\zp.\]
\noindent
That is,  $(a^{{t}_{a}}, a^{t_{b}})$ is a nontrivial type-1 solution of the equation (\ref{eq:BasicZ}) with $n=d$, since \\
 \[ x^d + 1 = y^d,  \text{ in }\zp. \] \\
Moreover,  if $e$ is any divisor of $d$; that is, $d = ee_{1}$, then 
$(a^{t_{a}e_{1}}, a^{t_{b}e_{1}})$ is a nontrivial solution type-1 of the equation (\ref{eq:BasicZ}) with $n=e$, since 
\[x^{e} + 1 = y^{e},  \text{ in }\zp. \]\\
Lastly, we note that we only need to consider $1\leq \alpha, \beta \leq \frac{(p-1)^2}{2}$. This follows because  we are only interested in the exponents $1 \leq  d \leq \frac{p-1}{2}$, while we can take $1\leq t_a, t_b \leq (p-1)$, because of the cyclic nature of $\zpstar$.

Here is a schematic presentation of one step of the algorithm. Notice that it involves $(p-1)^2/4$ applications of the Euclidean Algorithm.
We will have to repeat this for every $1 \leq \alpha \leq (p-1)$.

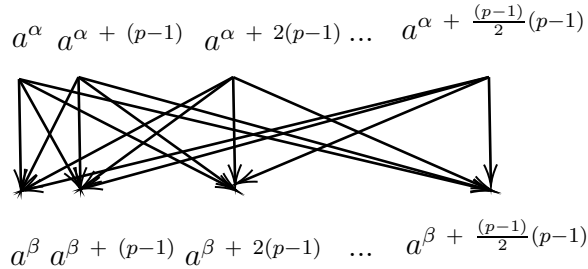
\tikzset{every picture/.style={line width=1pt}} 

\begin{tikzpicture}[x=0.75pt,y=0.75pt,yscale=-1,xscale=1]

\draw    (63.5,64) -- (64.46,118.5) ;
\draw [shift={(64.5,120.5)}, rotate = 268.99] [color={rgb, 255:red, 0; green, 0; blue, 0 }  ][line width=0.75]    (10.93,-3.29) .. controls (6.95,-1.4) and (3.31,-0.3) .. (0,0) .. controls (3.31,0.3) and (6.95,1.4) .. (10.93,3.29)   ;

\draw    (63.5,64) -- (93.52,117.75) ;
\draw [shift={(94.5,119.5)}, rotate = 240.81] [color={rgb, 255:red, 0; green, 0; blue, 0 }  ][line width=0.75]    (10.93,-3.29) .. controls (6.95,-1.4) and (3.31,-0.3) .. (0,0) .. controls (3.31,0.3) and (6.95,1.4) .. (10.93,3.29)   ;

\draw    (63.5,64) -- (171.72,119.59) ;
\draw [shift={(173.5,120.5)}, rotate = 207.19] [color={rgb, 255:red, 0; green, 0; blue, 0 }  ][line width=0.75]    (10.93,-3.29) .. controls (6.95,-1.4) and (3.31,-0.3) .. (0,0) .. controls (3.31,0.3) and (6.95,1.4) .. (10.93,3.29)   ;

\draw    (63.5,64) -- (300.55,120.04) ;
\draw [shift={(302.5,120.5)}, rotate = 193.3] [color={rgb, 255:red, 0; green, 0; blue, 0 }  ][line width=0.75]    (10.93,-3.29) .. controls (6.95,-1.4) and (3.31,-0.3) .. (0,0) .. controls (3.31,0.3) and (6.95,1.4) .. (10.93,3.29)   ;

\draw    (93.5,63) -- (94.46,117.5) ;
\draw [shift={(94.5,119.5)}, rotate = 268.99] [color={rgb, 255:red, 0; green, 0; blue, 0 }  ][line width=0.75]    (10.93,-3.29) .. controls (6.95,-1.4) and (3.31,-0.3) .. (0,0) .. controls (3.31,0.3) and (6.95,1.4) .. (10.93,3.29)   ;

\draw    (93.5,63) -- (171.88,119.33) ;
\draw [shift={(173.5,120.5)}, rotate = 215.71] [color={rgb, 255:red, 0; green, 0; blue, 0 }  ][line width=0.75]    (10.93,-3.29) .. controls (6.95,-1.4) and (3.31,-0.3) .. (0,0) .. controls (3.31,0.3) and (6.95,1.4) .. (10.93,3.29)   ;

\draw    (93.5,63) -- (65.4,118.71) ;
\draw [shift={(64.5,120.5)}, rotate = 296.76] [color={rgb, 255:red, 0; green, 0; blue, 0 }  ][line width=0.75]    (10.93,-3.29) .. controls (6.95,-1.4) and (3.31,-0.3) .. (0,0) .. controls (3.31,0.3) and (6.95,1.4) .. (10.93,3.29)   ;

\draw    (93.5,63) -- (300.57,119.97) ;
\draw [shift={(302.5,120.5)}, rotate = 195.38] [color={rgb, 255:red, 0; green, 0; blue, 0 }  ][line width=0.75]    (10.93,-3.29) .. controls (6.95,-1.4) and (3.31,-0.3) .. (0,0) .. controls (3.31,0.3) and (6.95,1.4) .. (10.93,3.29)   ;

\draw    (171.5,63) -- (172.46,116.03) ;
\draw [shift={(172.5,118.03)}, rotate = 268.96] [color={rgb, 255:red, 0; green, 0; blue, 0 }  ][line width=0.75]    (10.93,-3.29) .. controls (6.95,-1.4) and (3.31,-0.3) .. (0,0) .. controls (3.31,0.3) and (6.95,1.4) .. (10.93,3.29)   ;

\draw    (171.5,63) -- (96.11,118.32) ;
\draw [shift={(94.5,119.5)}, rotate = 323.73] [color={rgb, 255:red, 0; green, 0; blue, 0 }  ][line width=0.75]    (10.93,-3.29) .. controls (6.95,-1.4) and (3.31,-0.3) .. (0,0) .. controls (3.31,0.3) and (6.95,1.4) .. (10.93,3.29)   ;

\draw    (171.5,63) -- (66.26,119.55) ;
\draw [shift={(64.5,120.5)}, rotate = 331.75] [color={rgb, 255:red, 0; green, 0; blue, 0 }  ][line width=0.75]    (10.93,-3.29) .. controls (6.95,-1.4) and (3.31,-0.3) .. (0,0) .. controls (3.31,0.3) and (6.95,1.4) .. (10.93,3.29)   ;

\draw    (171.5,63) -- (300.67,119.7) ;
\draw [shift={(302.5,120.5)}, rotate = 203.7] [color={rgb, 255:red, 0; green, 0; blue, 0 }  ][line width=0.75]    (10.93,-3.29) .. controls (6.95,-1.4) and (3.31,-0.3) .. (0,0) .. controls (3.31,0.3) and (6.95,1.4) .. (10.93,3.29)   ;

\draw    (300.5,63) -- (301.46,116.03) ;
\draw [shift={(301.5,118.03)}, rotate = 268.96] [color={rgb, 255:red, 0; green, 0; blue, 0 }  ][line width=0.75]    (10.93,-3.29) .. controls (6.95,-1.4) and (3.31,-0.3) .. (0,0) .. controls (3.31,0.3) and (6.95,1.4) .. (10.93,3.29)   ;

\draw    (300.5,63) -- (96.43,118.97) ;
\draw [shift={(94.5,119.5)}, rotate = 344.65999999999997] [color={rgb, 255:red, 0; green, 0; blue, 0 }  ][line width=0.75]    (10.93,-3.29) .. controls (6.95,-1.4) and (3.31,-0.3) .. (0,0) .. controls (3.31,0.3) and (6.95,1.4) .. (10.93,3.29)   ;

\draw    (300.5,63) -- (174.34,116.25) ;
\draw [shift={(172.5,117.03)}, rotate = 337.12] [color={rgb, 255:red, 0; green, 0; blue, 0 }  ][line width=0.75]    (10.93,-3.29) .. controls (6.95,-1.4) and (3.31,-0.3) .. (0,0) .. controls (3.31,0.3) and (6.95,1.4) .. (10.93,3.29)   ;

\draw    (300.5,63) -- (66.44,120.03) ;
\draw [shift={(64.5,120.5)}, rotate = 346.31] [color={rgb, 255:red, 0; green, 0; blue, 0 }  ][line width=0.75]    (10.93,-3.29) .. controls (6.95,-1.4) and (3.31,-0.3) .. (0,0) .. controls (3.31,0.3) and (6.95,1.4) .. (10.93,3.29)   ;

\draw (721,21) node   {$0$};
\draw (701,71) node   {$0$};
\draw (68,42.5) node   {$a^{\alpha }$};
\draw (115,42.5) node   {$a^{\alpha \ +\ ( p-1)}$};
\draw (192,42.5) node   {$a^{\alpha \ +\ 2( p-1)}$};
\draw (303,35) node   {$a^{\alpha \ +\ \frac{( p-1)}{2}( p-1)}$};
\draw (236,44) node  [align=left] {...};
\draw (67,151) node   {$a^{\beta }$};
\draw (111,151) node   {$a^{\beta \ +\ ( p-1)}$};
\draw (182,151) node   {$a^{\beta \ +\ 2( p-1)}$};
\draw (305,143.5) node   {$a^{\beta \ +\ \frac{( p-1)}{2}( p-1)}$};
\draw (236,152.5) node  [align=left] {...};

\end{tikzpicture}

\section{Examples}\label{sec:E}
Below we show two examples for $p = 17 = 2^4 + 1$ and $p= 23 = 2*11 + 1$. The solutions were computed using a Python implementation of the algorithm. We chose these primes because the structure of $\mathbb{Z}^*_{17}$ and $\mathbb{Z}^*_{23}$ are essentially different in ways which are typical.  The layout of the tables should make clear the structural relationships in the solutions for different powers $n$ which are detailed in sections (\ref{sec:T0}) and (\ref{sec:T1}).

Recall that if $x^n + y^n = z^n$ then for any non-zero element $u \in \zpstar$,  we also have $(u\cdot x)^n + (u\cdot y)^n = (u \cdot z)^n$. Consequently we only need to find the basics type-0 solutions $x^n + 1 = 0$ and the basic type-1 solutions $x^n + 1 = z^n$.

Example 1: $p = 17, p-1 = 2^4.$\\
In each odd power $n$, we will only have one basic type-0 solution: 
$$16^n + 1 = 0.$$
\begin{table}[h]
\tiny
\caption{\label{tab:Type0_p17_n_even} Basic type-0 solutions in even powers $\mathbb{Z}_{17}$ .}
\begin{tabular}{| rl | rl | rl | rl |}
\hline
				n&=2	&			n&=2*3	&				n&=2*5	&				n&=2*7\\
\hline
$13^{2} 	+  1$ 	&$= 0$ 	&$13^{6} +  1$	&$= 0$ 	& $13^{10} +  1$  	&$= 0$ 	&  $13^{14} +  1 $ 	& $= 0$\\
$4^{2} 	+  1 $ 	&$= 0$  	&$4^{6} +  1$	&$= 0$  	& $4^{10} +  1 $ 	& $= 0$ 	& $4^{14} +  1$ 	& $= 0$\\
\hline
				n&=4	&			n&=4*3	&				&		&				&\\
\hline
$9^{4} 	+  1$  	&$= 0$   	&$9^{12} +  1$ 	&$= 0$	&				&		&				&\\
$8^{4} 	+  1$ 	&$= 0$	&$8^{12} + 1$ 	&$= 0$	&				&		&				&\\
\hline
				n&=8	&			&		&				&		&				&\\
\hline
$3^{8} 	+  1$		&$= 0$   	&			&		&				&		&				&\\
$14^{8} 	+  1$		&$ = 0$	&			&		&				&		&				&\\
\hline
\end{tabular}
\end{table}

As for the basic type-1 solutions, we list them below in two separate tables:   for even then for odd $n$'s.\\
\begin{table}[h]
\tiny
\caption{\label{tab:Type1_p17_n_even} Basic type-1 solutions in even powers $n$ in $\mathbb{Z}_{17}$ .}
\begin{tabular}{| rlrl | rlrl |}
\hline
 & &n = 2& &  & &n = 14& \\

\hline
    $1^{2} +  1 $ &$= 6^{2}$  &  $16^{2} +  1 $ &$= 11^{2}$  &  $1^{14} +  1 $ &$= 3^{14}$  &  $16^{14} +  1 $ &$= 14^{14}$ \\  
    $1^{2} +  1 $ &$= 11^{2}$  &  $16^{2} +  1 $ &$= 6^{2}$  &  $1^{14} +  1 $ &$= 14^{14}$  &  $16^{14} +  1 $ &$= 3^{14}$ \\  
    $5^{2} +  1 $ &$= 3^{2}$  &  $12^{2} +  1 $ &$= 14^{2}$  &  $7^{14} +  1 $ &$= 6^{14}$  &  $10^{14} +  1 $ &$= 11^{14}$ \\  
    $5^{2} +  1 $ &$= 14^{2}$  &  $12^{2} +  1 $ &$= 3^{2}$  &  $7^{14} +  1 $ &$= 11^{14}$  &  $10^{14} +  1 $ &$= 6^{14}$ \\  
    $7^{2} +  1 $ &$= 4^{2}$  &  $10^{2} +  1 $ &$= 13^{2}$  &  $5^{14} +  1 $ &$= 13^{14}$  &  $12^{14} +  1 $ &$= 4^{14}$ \\  
    $7^{2} +  1 $ &$= 13^{2}$  &  $10^{2} +  1 $ &$= 4^{2}$  &  $5^{14} +  1 $ &$= 4^{14}$  &  $12^{14} +  1 $ &$= 13^{14}$ \\  
\hline
 & &n = 6& &  & &n = 10& \\

\hline
    $1^{6} +  1 $ &$= 5^{6}$  &  $16^{6} +  1 $ &$= 12^{6}$  &  $1^{10} +  1 $ &$= 7^{10}$  &  $16^{10} +  1 $ &$= 10^{10}$ \\  
    $1^{6} +  1 $ &$= 12^{6}$  &  $16^{6} +  1 $ &$= 5^{6}$  &  $1^{10} +  1 $ &$= 10^{10}$  &  $16^{10} +  1 $ &$= 7^{10}$ \\  
    $3^{6} +  1 $ &$= 4^{6}$  &  $14^{6} +  1 $ &$= 13^{6}$  &  $6^{10} +  1 $ &$= 13^{10}$  &  $11^{10} +  1 $ &$= 4^{10}$ \\  
    $3^{6} +  1 $ &$= 13^{6}$  &  $14^{6} +  1 $ &$= 4^{6}$  &  $6^{10} +  1 $ &$= 4^{10}$  &  $11^{10} +  1 $ &$= 13^{10}$ \\  
    $6^{6} +  1 $ &$= 7^{6}$  &  $11^{6} +  1 $ &$= 10^{6}$  &  $3^{10} +  1 $ &$= 5^{10}$  &  $14^{10} +  1 $ &$= 12^{10}$ \\  
    $6^{6} +  1 $ &$= 10^{6}$  &  $11^{6} +  1 $ &$= 7^{6}$  &  $3^{10} +  1 $ &$= 12^{10}$  &  $14^{10} +  1 $ &$= 5^{10}$ \\ 
    \hline
\end{tabular}

\end{table}
Notice that there are no type-1 solutions in powers 4, 8, 12, 16.\\

\newpage
\begin{table}[h]
\tiny
\caption{\label{tab:Type1_p17_n_odd}Basic type-1 solutions in odd powers $n$ in $\mathbb{Z}_{17}$ .}
\begin{tabular}{| rlrl | rlrl |}
\hline
 & &n = 1& &  & &n = 15& \\
 \hline
    $1^{1} +  1 $ &$= 2^{1}$  &  $15^{1} +  1 $ &$= 16^{1}$  &  $1^{15} +  1 $ &$= 9^{15}$  &  $8^{15} +  1 $ &$= 16^{15}$ \\  
    $2^{1} +  1 $ &$= 3^{1}$  &  $14^{1} +  1 $ &$= 15^{1}$  &  $9^{15} +  1 $ &$= 6^{15}$  &  $11^{15} +  1 $ &$= 8^{15}$ \\  
    $3^{1} +  1 $ &$= 4^{1}$  &  $13^{1} +  1 $ &$= 14^{1}$  &  $6^{15} +  1 $ &$= 13^{15}$  &  $4^{15} +  1 $ &$= 11^{15}$ \\  
    $4^{1} +  1 $ &$= 5^{1}$  &  $12^{1} +  1 $ &$= 13^{1}$  &  $13^{15} +  1 $ &$= 7^{15}$  &  $10^{15} +  1 $ &$= 4^{15}$ \\  
    $5^{1} +  1 $ &$= 6^{1}$  &  $11^{1} +  1 $ &$= 12^{1}$  &  $7^{15} +  1 $ &$= 3^{15}$  &  $14^{15} +  1 $ &$= 10^{15}$ \\  
    $6^{1} +  1 $ &$= 7^{1}$  &  $10^{1} +  1 $ &$= 11^{1}$  &  $3^{15} +  1 $ &$= 5^{15}$  &  $12^{15} +  1 $ &$= 14^{15}$ \\  
    $7^{1} +  1 $ &$= 8^{1}$  &  $9^{1} +  1 $ &$= 10^{1}$  &  $5^{15} +  1 $ &$= 15^{15}$  &  $2^{15} +  1 $ &$= 12^{15}$ \\  
    $8^{1} +  1 $ &$= 9^{1}$  &  $8^{1} +  1 $ &$= 9^{1}$  &  $15^{15} +  1 $ &$= 2^{15}$  &  $15^{15} +  1 $ &$= 2^{15}$ \\ 
\hline
 & &n = 3& &  & &n = 13& \\
\hline
    $1^{3} +  1 $ &$= 8^{3}$  &  $9^{3} +  1 $ &$= 16^{3}$  &  $1^{13} +  1 $ &$= 15^{13}$  &  $2^{13} +  1 $ &$= 16^{13}$ \\  
    $2^{3} +  1 $ &$= 15^{3}$  &  $2^{3} +  1 $ &$= 15^{3}$  &  $9^{13} +  1 $ &$= 8^{13}$  &  $9^{13} +  1 $ &$= 8^{13}$ \\  
    $3^{3} +  1 $ &$= 12^{3}$  &  $5^{3} +  1 $ &$= 14^{3}$  &  $6^{13} +  1 $ &$= 10^{13}$  &  $7^{13} +  1 $ &$= 11^{13}$ \\  
    $4^{3} +  1 $ &$= 10^{3}$  &  $7^{3} +  1 $ &$= 13^{3}$  &  $13^{13} +  1 $ &$= 12^{13}$  &  $5^{13} +  1 $ &$= 4^{13}$ \\  
    $5^{3} +  1 $ &$= 14^{3}$  &  $3^{3} +  1 $ &$= 12^{3}$  &  $7^{13} +  1 $ &$= 11^{13}$  &  $6^{13} +  1 $ &$= 10^{13}$ \\  
    $6^{3} +  1 $ &$= 4^{3}$  &  $13^{3} +  1 $ &$= 11^{3}$  &  $3^{13} +  1 $ &$= 13^{13}$  &  $4^{13} +  1 $ &$= 14^{13}$ \\  
    $7^{3} +  1 $ &$= 13^{3}$  &  $4^{3} +  1 $ &$= 10^{3}$  &  $5^{13} +  1 $ &$= 4^{13}$  &  $13^{13} +  1 $ &$= 12^{13}$ \\  
    $8^{3} +  1 $ &$= 7^{3}$  &  $10^{3} +  1 $ &$= 9^{3}$  &  $15^{13} +  1 $ &$= 5^{13}$  &  $12^{13} +  1 $ &$= 2^{13}$ \\ 
\hline
     & &n = 5& &  & &n = 11& \\
\hline
    $1^{5} +  1 $ &$= 15^{5}$  &  $2^{5} +  1 $ &$= 16^{5}$  &  $1^{11} +  1 $ &$= 8^{11}$  &  $9^{11} +  1 $ &$= 16^{11}$ \\  
    $2^{5} +  1 $ &$= 16^{5}$  &  $1^{5} +  1 $ &$= 15^{5}$  &  $9^{11} +  1 $ &$= 16^{11}$  &  $1^{11} +  1 $ &$= 8^{11}$ \\  
    $3^{5} +  1 $ &$= 10^{5}$  &  $7^{5} +  1 $ &$= 14^{5}$  &  $6^{11} +  1 $ &$= 12^{11}$  &  $5^{11} +  1 $ &$= 11^{11}$ \\  
    $4^{5} +  1 $ &$= 3^{5}$  &  $14^{5} +  1 $ &$= 13^{5}$  &  $13^{11} +  1 $ &$= 6^{11}$  &  $11^{11} +  1 $ &$= 4^{11}$ \\  
    $5^{5} +  1 $ &$= 2^{5}$  &  $15^{5} +  1 $ &$= 12^{5}$  &  $7^{11} +  1 $ &$= 9^{11}$  &  $8^{11} +  1 $ &$= 10^{11}$ \\  
    $6^{5} +  1 $ &$= 9^{5}$  &  $8^{5} +  1 $ &$= 11^{5}$  &  $3^{11} +  1 $ &$= 2^{11}$  &  $15^{11} +  1 $ &$= 14^{11}$ \\  
    $7^{5} +  1 $ &$= 14^{5}$  &  $3^{5} +  1 $ &$= 10^{5}$  &  $5^{11} +  1 $ &$= 11^{11}$  &  $6^{11} +  1 $ &$= 12^{11}$ \\  
    $8^{5} +  1 $ &$= 11^{5}$  &  $6^{5} +  1 $ &$= 9^{5}$  &  $15^{11} +  1 $ &$= 14^{11}$  &  $3^{11} +  1 $ &$= 2^{11}$ \\ 
 \hline   
     & &n = 7& &  & &n = 9& \\
\hline
    $1^{7} +  1 $ &$= 9^{7}$  &  $8^{7} +  1 $ &$= 16^{7}$  &  $1^{9} +  1 $ &$= 2^{9}$  &  $15^{9} +  1 $ &$= 16^{9}$ \\  
    $2^{7} +  1 $ &$= 5^{7}$  &  $12^{7} +  1 $ &$= 15^{7}$  &  $9^{9} +  1 $ &$= 7^{9}$  &  $10^{9} +  1 $ &$= 8^{9}$ \\  
    $3^{7} +  1 $ &$= 7^{7}$  &  $10^{7} +  1 $ &$= 14^{7}$  &  $6^{9} +  1 $ &$= 5^{9}$  &  $12^{9} +  1 $ &$= 11^{9}$ \\  
    $4^{7} +  1 $ &$= 6^{7}$  &  $11^{7} +  1 $ &$= 13^{7}$  &  $13^{9} +  1 $ &$= 3^{9}$  &  $14^{9} +  1 $ &$= 4^{9}$ \\  
    $5^{7} +  1 $ &$= 3^{7}$  &  $14^{7} +  1 $ &$= 12^{7}$  &  $7^{9} +  1 $ &$= 6^{9}$  &  $11^{9} +  1 $ &$= 10^{9}$ \\  
    $6^{7} +  1 $ &$= 8^{7}$  &  $9^{7} +  1 $ &$= 11^{7}$  &  $3^{9} +  1 $ &$= 15^{9}$  &  $2^{9} +  1 $ &$= 14^{9}$ \\  
    $7^{7} +  1 $ &$= 4^{7}$  &  $13^{7} +  1 $ &$= 10^{7}$  &  $5^{9} +  1 $ &$= 13^{9}$  &  $4^{9} +  1 $ &$= 12^{9}$ \\  
    $8^{7} +  1 $ &$= 16^{7}$  &  $1^{7} +  1 $ &$= 9^{7}$  &  $15^{9} +  1 $ &$= 16^{9}$  &  $1^{9} +  1 $ &$= 2^{9}$ \\
 \hline  
\end{tabular}

\end{table}

Example 2:  $p = 23, p-1 = 2*11.$\\
We have no type-0 solutions in any even power and, in each odd power $n$, we only have one basic type-0 solution:
 $$22^n + 1 = 0.$$\\
For type-1 solutions, we list all the basic solutions, first in even, then in odd powers $n$:\\
\begin{table}[h]
\tiny
\caption{\label{tab:Type1_p23_n_even}Basic type-1 solutions in even powers $n$ in $\mathbb{Z}_{23}$ .}
\begin{tabular}{| rlrl | rlrl |}
\hline
 & &n = 2& &  & &n = 20& \\
\hline
    $1^{2} +  1 $ &$= 5^{2}$  &  $22^{2} +  1 $ &$= 18^{2}$  &  $1^{20} +  1 $ &$= 14^{20}$  &  $22^{20} +  1 $ &$= 9^{20}$ \\  
    $1^{2} +  1 $ &$= 18^{2}$  &  $22^{2} +  1 $ &$= 5^{2}$  &  $1^{20} +  1 $ &$= 9^{20}$  &  $22^{20} +  1 $ &$= 14^{20}$ \\  
    $5^{2} +  1 $ &$= 7^{2}$  &  $18^{2} +  1 $ &$= 16^{2}$  &  $14^{20} +  1 $ &$= 10^{20}$  &  $9^{20} +  1 $ &$= 13^{20}$ \\  
    $5^{2} +  1 $ &$= 16^{2}$  &  $18^{2} +  1 $ &$= 7^{2}$  &  $14^{20} +  1 $ &$= 13^{20}$  &  $9^{20} +  1 $ &$= 10^{20}$ \\  
    $7^{2} +  1 $ &$= 2^{2}$  &  $16^{2} +  1 $ &$= 21^{2}$  &  $10^{20} +  1 $ &$= 12^{20}$  &  $13^{20} +  1 $ &$= 11^{20}$ \\  
    $7^{2} +  1 $ &$= 21^{2}$  &  $16^{2} +  1 $ &$= 2^{2}$  &  $10^{20} +  1 $ &$= 11^{20}$  &  $13^{20} +  1 $ &$= 12^{20}$ \\  
    $9^{2} +  1 $ &$= 6^{2}$  &  $14^{2} +  1 $ &$= 17^{2}$  &  $18^{20} +  1 $ &$= 4^{20}$  &  $5^{20} +  1 $ &$= 19^{20}$ \\  
    $9^{2} +  1 $ &$= 17^{2}$  &  $14^{2} +  1 $ &$= 6^{2}$  &  $18^{20} +  1 $ &$= 19^{20}$  &  $5^{20} +  1 $ &$= 4^{20}$ \\  
    $10^{2} +  1 $ &$= 3^{2}$  &  $13^{2} +  1 $ &$= 20^{2}$  &  $7^{20} +  1 $ &$= 8^{20}$  &  $16^{20} +  1 $ &$= 15^{20}$ \\  
    $10^{2} +  1 $ &$= 20^{2}$  &  $13^{2} +  1 $ &$= 3^{2}$  &  $7^{20} +  1 $ &$= 15^{20}$  &  $16^{20} +  1 $ &$= 8^{20}$ \\
\hline
& &n = 4& &  & &n = 18& \\

\hline
    $1^{4} +  1 $ &$= 8^{4}$  &  $22^{4} +  1 $ &$= 15^{4}$  &  $1^{18} +  1 $ &$= 3^{18}$  &  $22^{18} +  1 $ &$= 20^{18}$ \\  
    $1^{4} +  1 $ &$= 15^{4}$  &  $22^{4} +  1 $ &$= 8^{4}$  &  $1^{18} +  1 $ &$= 20^{18}$  &  $22^{18} +  1 $ &$= 3^{18}$ \\  
    $3^{4} +  1 $ &$= 11^{4}$  &  $20^{4} +  1 $ &$= 12^{4}$  &  $8^{18} +  1 $ &$= 21^{18}$  &  $15^{18} +  1 $ &$= 2^{18}$ \\  
    $3^{4} +  1 $ &$= 12^{4}$  &  $20^{4} +  1 $ &$= 11^{4}$  &  $8^{18} +  1 $ &$= 2^{18}$  &  $15^{18} +  1 $ &$= 21^{18}$ \\  
    $4^{4} +  1 $ &$= 5^{4}$  &  $19^{4} +  1 $ &$= 18^{4}$  &  $6^{18} +  1 $ &$= 14^{18}$  &  $17^{18} +  1 $ &$= 9^{18}$ \\  
    $4^{4} +  1 $ &$= 18^{4}$  &  $19^{4} +  1 $ &$= 5^{4}$  &  $6^{18} +  1 $ &$= 9^{18}$  &  $17^{18} +  1 $ &$= 14^{18}$ \\  
    $6^{4} +  1 $ &$= 7^{4}$  &  $17^{4} +  1 $ &$= 16^{4}$  &  $4^{18} +  1 $ &$= 10^{18}$  &  $19^{18} +  1 $ &$= 13^{18}$ \\  
    $6^{4} +  1 $ &$= 16^{4}$  &  $17^{4} +  1 $ &$= 7^{4}$  &  $4^{18} +  1 $ &$= 13^{18}$  &  $19^{18} +  1 $ &$= 10^{18}$ \\  
    $8^{4} +  1 $ &$= 4^{4}$  &  $15^{4} +  1 $ &$= 19^{4}$  &  $3^{18} +  1 $ &$= 6^{18}$  &  $20^{18} +  1 $ &$= 17^{18}$ \\  
    $8^{4} +  1 $ &$= 19^{4}$  &  $15^{4} +  1 $ &$= 4^{4}$  &  $3^{18} +  1 $ &$= 17^{18}$  &  $20^{18} +  1 $ &$= 6^{18}$ \\
\hline
 & &n = 6& &  & &n = 16& \\

\hline
    $1^{6} +  1 $ &$= 4^{6}$  &  $22^{6} +  1 $ &$= 19^{6}$  &  $1^{16} +  1 $ &$= 6^{16}$  &  $22^{16} +  1 $ &$= 17^{16}$ \\  
    $1^{6} +  1 $ &$= 19^{6}$  &  $22^{6} +  1 $ &$= 4^{6}$  &  $1^{16} +  1 $ &$= 17^{16}$  &  $22^{16} +  1 $ &$= 6^{16}$ \\  
    $4^{6} +  1 $ &$= 9^{6}$  &  $19^{6} +  1 $ &$= 14^{6}$  &  $6^{16} +  1 $ &$= 18^{16}$  &  $17^{16} +  1 $ &$= 5^{16}$ \\  
    $4^{6} +  1 $ &$= 14^{6}$  &  $19^{6} +  1 $ &$= 9^{6}$  &  $6^{16} +  1 $ &$= 5^{16}$  &  $17^{16} +  1 $ &$= 18^{16}$ \\  
    $5^{6} +  1 $ &$= 11^{6}$  &  $18^{6} +  1 $ &$= 12^{6}$  &  $14^{16} +  1 $ &$= 21^{16}$  &  $9^{16} +  1 $ &$= 2^{16}$ \\  
    $5^{6} +  1 $ &$= 12^{6}$  &  $18^{6} +  1 $ &$= 11^{6}$  &  $14^{16} +  1 $ &$= 2^{16}$  &  $9^{16} +  1 $ &$= 21^{16}$ \\  
    $6^{6} +  1 $ &$= 8^{6}$  &  $17^{6} +  1 $ &$= 15^{6}$  &  $4^{16} +  1 $ &$= 3^{16}$  &  $19^{16} +  1 $ &$= 20^{16}$ \\  
    $6^{6} +  1 $ &$= 15^{6}$  &  $17^{6} +  1 $ &$= 8^{6}$  &  $4^{16} +  1 $ &$= 20^{16}$  &  $19^{16} +  1 $ &$= 3^{16}$ \\  
    $9^{6} +  1 $ &$= 7^{6}$  &  $14^{6} +  1 $ &$= 16^{6}$  &  $18^{16} +  1 $ &$= 10^{16}$  &  $5^{16} +  1 $ &$= 13^{16}$ \\  
    $9^{6} +  1 $ &$= 16^{6}$  &  $14^{6} +  1 $ &$= 7^{6}$  &  $18^{16} +  1 $ &$= 13^{16}$  &  $5^{16} +  1 $ &$= 10^{16}$ \\
 \hline
  & &n = 8& &  & &n = 14& \\

\hline
    $1^{8} +  1 $ &$= 10^{8}$  &  $22^{8} +  1 $ &$= 13^{8}$  &  $1^{14} +  1 $ &$= 7^{14}$  &  $22^{14} +  1 $ &$= 16^{14}$ \\  
    $1^{8} +  1 $ &$= 13^{8}$  &  $22^{8} +  1 $ &$= 10^{8}$  &  $1^{14} +  1 $ &$= 16^{14}$  &  $22^{14} +  1 $ &$= 7^{14}$ \\  
    $2^{8} +  1 $ &$= 8^{8}$  &  $21^{8} +  1 $ &$= 15^{8}$  &  $12^{14} +  1 $ &$= 3^{14}$  &  $11^{14} +  1 $ &$= 20^{14}$ \\  
    $2^{8} +  1 $ &$= 15^{8}$  &  $21^{8} +  1 $ &$= 8^{8}$  &  $12^{14} +  1 $ &$= 20^{14}$  &  $11^{14} +  1 $ &$= 3^{14}$ \\  
    $7^{8} +  1 $ &$= 9^{8}$  &  $16^{8} +  1 $ &$= 14^{8}$  &  $10^{14} +  1 $ &$= 18^{14}$  &  $13^{14} +  1 $ &$= 5^{14}$ \\  
    $7^{8} +  1 $ &$= 14^{8}$  &  $16^{8} +  1 $ &$= 9^{8}$  &  $10^{14} +  1 $ &$= 5^{14}$  &  $13^{14} +  1 $ &$= 18^{14}$ \\  
    $10^{8} +  1 $ &$= 2^{8}$  &  $13^{8} +  1 $ &$= 21^{8}$  &  $7^{14} +  1 $ &$= 12^{14}$  &  $16^{14} +  1 $ &$= 11^{14}$ \\  
    $10^{8} +  1 $ &$= 21^{8}$  &  $13^{8} +  1 $ &$= 2^{8}$  &  $7^{14} +  1 $ &$= 11^{14}$  &  $16^{14} +  1 $ &$= 12^{14}$ \\  
    $11^{8} +  1 $ &$= 4^{8}$  &  $12^{8} +  1 $ &$= 19^{8}$  &  $21^{14} +  1 $ &$= 6^{14}$  &  $2^{14} +  1 $ &$= 17^{14}$ \\  
    $11^{8} +  1 $ &$= 19^{8}$  &  $12^{8} +  1 $ &$= 4^{8}$  &  $21^{14} +  1 $ &$= 17^{14}$  &  $2^{14} +  1 $ &$= 6^{14}$ \\ 
 \hline
 & &n = 10& &  & &n = 12& \\
\hline
    $1^{10} +  1 $ &$= 11^{10}$  &  $22^{10} +  1 $ &$= 12^{10}$  &  $1^{12} +  1 $ &$= 21^{12}$  &  $22^{12} +  1 $ &$= 2^{12}$ \\  
    $1^{10} +  1 $ &$= 12^{10}$  &  $22^{10} +  1 $ &$= 11^{10}$  &  $1^{12} +  1 $ &$= 2^{12}$  &  $22^{12} +  1 $ &$= 21^{12}$ \\  
    $2^{10} +  1 $ &$= 7^{10}$  &  $21^{10} +  1 $ &$= 16^{10}$  &  $12^{12} +  1 $ &$= 10^{12}$  &  $11^{12} +  1 $ &$= 13^{12}$ \\  
    $2^{10} +  1 $ &$= 16^{10}$  &  $21^{10} +  1 $ &$= 7^{10}$  &  $12^{12} +  1 $ &$= 13^{12}$  &  $11^{12} +  1 $ &$= 10^{12}$ \\  
    $3^{10} +  1 $ &$= 5^{10}$  &  $20^{10} +  1 $ &$= 18^{10}$  &  $8^{12} +  1 $ &$= 14^{12}$  &  $15^{12} +  1 $ &$= 9^{12}$ \\  
    $3^{10} +  1 $ &$= 18^{10}$  &  $20^{10} +  1 $ &$= 5^{10}$  &  $8^{12} +  1 $ &$= 9^{12}$  &  $15^{12} +  1 $ &$= 14^{12}$ \\  
    $8^{10} +  1 $ &$= 6^{10}$  &  $15^{10} +  1 $ &$= 17^{10}$  &  $3^{12} +  1 $ &$= 4^{12}$  &  $20^{12} +  1 $ &$= 19^{12}$ \\  
    $8^{10} +  1 $ &$= 17^{10}$  &  $15^{10} +  1 $ &$= 6^{10}$  &  $3^{12} +  1 $ &$= 19^{12}$  &  $20^{12} +  1 $ &$= 4^{12}$ \\  
    $11^{10} +  1 $ &$= 8^{10}$  &  $12^{10} +  1 $ &$= 15^{10}$  &  $21^{12} +  1 $ &$= 3^{12}$  &  $2^{12} +  1 $ &$= 20^{12}$ \\  
    $11^{10} +  1 $ &$= 15^{10}$  &  $12^{10} +  1 $ &$= 8^{10}$  &  $21^{12} +  1 $ &$= 20^{12}$  &  $2^{12} +  1 $ &$= 3^{12}$ \\  
 \hline
\end{tabular}
\end{table}

\begin{table}[h]
\tiny
\caption{\label{tab:Type1_p23_n_odd}Basic type-1 solutions in odd powers $n$ in $\mathbb{Z}_{23}$ .}
\begin{tabular}{| rlrl | rlrl |}
\hline
&&n=1&    &   &&n=21&\\ 
\hline
 $1^{1} +  1$ &$= 2^{1}$ &$21^{1} +  1 $&$= 22^{1}$  & $1^{21} +  1 $ & $=12^{21}$   & $11^{21} +  1 $&$=22^{21}$ \\  
 $2^{1} +  1$ &$= 3^{1}$ &$20^{1} +  1 $&$= 21^{1}$  & $12^{21} +  1 $ & $=8^{21}$  &  $15^{21} +  1 $&$=11^{21}$ \\  
 $3^{1} +  1$ &$= 4^{1}$ &$19^{1} +  1 $&$= 20^{1}$  & $8^{21} +  1 $ & $=6^{21}$   &   $17^{21} +  1 $&$=15^{21}$ \\  
 $4^{1} +  1$ &$= 5^{1}$ &$18^{1} +  1 $&$= 19^{1}$  & $6^{21} +  1$ &  $= 14^{21}$  & $9^{21} +  1 $&$=17^{21}$ \\  
 $5^{1} +  1$ &$= 6^{1}$ &$17^{1} +  1 $&$= 18^{1}$  & $14^{21} +  1$ & $=4^{21}$  &  $19^{21} +  1 $&$=9^{21}$ \\  
 $6^{1} +  1$ &$= 7^{1}$ &$16^{1} +  1 $&$= 17^{1}$  & $4^{21} +  1 $ & $=10^{21}$  &  $13^{21} +  1 $&$=19^{21}$ \\  
 $7^{1} +  1$ &$= 8^{1}$ &$15^{1} +  1 $&$= 16^{1}$  & $10^{21} +  1 $ & $=3^{21}$  &  $20^{21} +  1 $&$=13^{21}$ \\  
 $8^{1} +  1$ &$= 9^{1}$ &$14^{1} +  1 $&$= 15^{1}$  & $3^{21} +  1 $ & $=18^{21}$   & $5^{21} +  1 $&$=20^{21}$ \\  
 $9^{1} +  1$ &$= 10^{1}$ &$13^{1} +  1 $&$= 14^{1}$  & $18^{21} +  1 $ & $=7^{21}$  &  $16^{21} +  1 $&$=5^{21}$ \\  
 $10^{1} +  1$ &$= 11^{1}$ &$12^{1} +  1 $&$= 13^{1}$  & $7^{21} +  1 $ & $=21^{21}$  &  $2^{21} +  1 $&$=16^{21}$ \\  
 $11^{1} +  1$ &$= 12^{1}$ &$11^{1} +  1 $&$= 12^{1}$  & $21^{21} +  1 $ & $=2^{21}$  &  $21^{21} +  1 $&$=2^{21}$ \\ 
\hline
& &n = 3& &  & &n = 19& \\
\hline
    $1^{3} +  1 $ &$= 16^{3}$  &  $7^{3} +  1 $ &$= 22^{3}$  &  $1^{19} +  1 $ &$= 13^{19}$  &  $10^{19} +  1 $ &$= 22^{19}$ \\  
    $2^{3} +  1 $ &$= 6^{3}$  &  $17^{3} +  1 $ &$= 21^{3}$  &  $12^{19} +  1 $ &$= 4^{19}$  &  $19^{19} +  1 $ &$= 11^{19}$ \\  
    $3^{3} +  1 $ &$= 19^{3}$  &  $4^{3} +  1 $ &$= 20^{3}$  &  $8^{19} +  1 $ &$= 17^{19}$  &  $6^{19} +  1 $ &$= 15^{19}$ \\  
    $4^{3} +  1 $ &$= 20^{3}$  &  $3^{3} +  1 $ &$= 19^{3}$  &  $6^{19} +  1 $ &$= 15^{19}$  &  $8^{19} +  1 $ &$= 17^{19}$ \\  
    $5^{3} +  1 $ &$= 10^{3}$  &  $13^{3} +  1 $ &$= 18^{3}$  &  $14^{19} +  1 $ &$= 7^{19}$  &  $16^{19} +  1 $ &$= 9^{19}$ \\  
    $6^{3} +  1 $ &$= 5^{3}$  &  $18^{3} +  1 $ &$= 17^{3}$  &  $4^{19} +  1 $ &$= 14^{19}$  &  $9^{19} +  1 $ &$= 19^{19}$ \\  
    $7^{3} +  1 $ &$= 22^{3}$  &  $1^{3} +  1 $ &$= 16^{3}$  &  $10^{19} +  1 $ &$= 22^{19}$  &  $1^{19} +  1 $ &$= 13^{19}$ \\  
    $8^{3} +  1 $ &$= 14^{3}$  &  $9^{3} +  1 $ &$= 15^{3}$  &  $3^{19} +  1 $ &$= 5^{19}$  &  $18^{19} +  1 $ &$= 20^{19}$ \\  
    $9^{3} +  1 $ &$= 15^{3}$  &  $8^{3} +  1 $ &$= 14^{3}$  &  $18^{19} +  1 $ &$= 20^{19}$  &  $3^{19} +  1 $ &$= 5^{19}$ \\  
    $10^{3} +  1 $ &$= 13^{3}$  &  $10^{3} +  1 $ &$= 13^{3}$  &  $7^{19} +  1 $ &$= 16^{19}$  &  $7^{19} +  1 $ &$= 16^{19}$ \\  
    $11^{3} +  1 $ &$= 7^{3}$  &  $16^{3} +  1 $ &$= 12^{3}$  &  $21^{19} +  1 $ &$= 10^{19}$  &  $13^{19} +  1 $ &$= 2^{19}$ \\ 
\hline    
     & &n = 5& &  & &n = 17& \\
\hline
    $1^{5} +  1 $ &$= 6^{5}$  &  $17^{5} +  1 $ &$= 22^{5}$  &  $1^{17} +  1 $ &$= 4^{17}$  &  $19^{17} +  1 $ &$= 22^{17}$ \\  
    $2^{5} +  1 $ &$= 20^{5}$  &  $3^{5} +  1 $ &$= 21^{5}$  &  $12^{17} +  1 $ &$= 15^{17}$  &  $8^{17} +  1 $ &$= 11^{17}$ \\  
    $3^{5} +  1 $ &$= 21^{5}$  &  $2^{5} +  1 $ &$= 20^{5}$  &  $8^{17} +  1 $ &$= 11^{17}$  &  $12^{17} +  1 $ &$= 15^{17}$ \\  
    $4^{5} +  1 $ &$= 3^{5}$  &  $20^{5} +  1 $ &$= 19^{5}$  &  $6^{17} +  1 $ &$= 8^{17}$  &  $15^{17} +  1 $ &$= 17^{17}$ \\  
    $5^{5} +  1 $ &$= 17^{5}$  &  $6^{5} +  1 $ &$= 18^{5}$  &  $14^{17} +  1 $ &$= 19^{17}$  &  $4^{17} +  1 $ &$= 9^{17}$ \\  
    $6^{5} +  1 $ &$= 18^{5}$  &  $5^{5} +  1 $ &$= 17^{5}$  &  $4^{17} +  1 $ &$= 9^{17}$  &  $14^{17} +  1 $ &$= 19^{17}$ \\  
    $7^{5} +  1 $ &$= 12^{5}$  &  $11^{5} +  1 $ &$= 16^{5}$  &  $10^{17} +  1 $ &$= 2^{17}$  &  $21^{17} +  1 $ &$= 13^{17}$ \\  
    $8^{5} +  1 $ &$= 7^{5}$  &  $16^{5} +  1 $ &$= 15^{5}$  &  $3^{17} +  1 $ &$= 10^{17}$  &  $13^{17} +  1 $ &$= 20^{17}$ \\  
    $9^{5} +  1 $ &$= 2^{5}$  &  $21^{5} +  1 $ &$= 14^{5}$  &  $18^{17} +  1 $ &$= 12^{17}$  &  $11^{17} +  1 $ &$= 5^{17}$ \\  
    $10^{5} +  1 $ &$= 5^{5}$  &  $18^{5} +  1 $ &$= 13^{5}$  &  $7^{17} +  1 $ &$= 14^{17}$  &  $9^{17} +  1 $ &$= 16^{17}$ \\  
    $11^{5} +  1 $ &$= 16^{5}$  &  $7^{5} +  1 $ &$= 12^{5}$  &  $21^{17} +  1 $ &$= 13^{17}$  &  $10^{17} +  1 $ &$= 2^{17}$ \\ 
 \hline
   & &n = 7& &  & &n = 15& \\
\hline
    $1^{7} +  1 $ &$= 3^{7}$  &  $20^{7} +  1 $ &$= 22^{7}$  &  $1^{15} +  1 $ &$= 8^{15}$  &  $15^{15} +  1 $ &$= 22^{15}$ \\  
    $2^{7} +  1 $ &$= 10^{7}$  &  $13^{7} +  1 $ &$= 21^{7}$  &  $12^{15} +  1 $ &$= 7^{15}$  &  $16^{15} +  1 $ &$= 11^{15}$ \\  
    $3^{7} +  1 $ &$= 6^{7}$  &  $17^{7} +  1 $ &$= 20^{7}$  &  $8^{15} +  1 $ &$= 4^{15}$  &  $19^{15} +  1 $ &$= 15^{15}$ \\  
    $4^{7} +  1 $ &$= 13^{7}$  &  $10^{7} +  1 $ &$= 19^{7}$  &  $6^{15} +  1 $ &$= 16^{15}$  &  $7^{15} +  1 $ &$= 17^{15}$ \\  
    $5^{7} +  1 $ &$= 16^{7}$  &  $7^{7} +  1 $ &$= 18^{7}$  &  $14^{15} +  1 $ &$= 13^{15}$  &  $10^{15} +  1 $ &$= 9^{15}$ \\  
    $6^{7} +  1 $ &$= 9^{7}$  &  $14^{7} +  1 $ &$= 17^{7}$  &  $4^{15} +  1 $ &$= 18^{15}$  &  $5^{15} +  1 $ &$= 19^{15}$ \\  
    $7^{7} +  1 $ &$= 18^{7}$  &  $5^{7} +  1 $ &$= 16^{7}$  &  $10^{15} +  1 $ &$= 9^{15}$  &  $14^{15} +  1 $ &$= 13^{15}$ \\  
    $8^{7} +  1 $ &$= 2^{7}$  &  $21^{7} +  1 $ &$= 15^{7}$  &  $3^{15} +  1 $ &$= 12^{15}$  &  $11^{15} +  1 $ &$= 20^{15}$ \\  
    $9^{7} +  1 $ &$= 7^{7}$  &  $16^{7} +  1 $ &$= 14^{7}$  &  $18^{15} +  1 $ &$= 10^{15}$  &  $13^{15} +  1 $ &$= 5^{15}$ \\  
    $10^{7} +  1 $ &$= 19^{7}$  &  $4^{7} +  1 $ &$= 13^{7}$  &  $7^{15} +  1 $ &$= 17^{15}$  &  $6^{15} +  1 $ &$= 16^{15}$ \\  
    $11^{7} +  1 $ &$= 4^{7}$  &  $19^{7} +  1 $ &$= 12^{7}$  &  $21^{15} +  1 $ &$= 6^{15}$  &  $17^{15} +  1 $ &$= 2^{15}$ \\ 
 \hline
& &n = 9& &  & &n = 13& \\
\hline
    $1^{9} +  1 $ &$= 9^{9}$  &  $14^{9} +  1 $ &$= 22^{9}$  &  $1^{13} +  1 $ &$= 18^{13}$  &  $5^{13} +  1 $ &$= 22^{13}$ \\  
    $2^{9} +  1 $ &$= 17^{9}$  &  $6^{9} +  1 $ &$= 21^{9}$  &  $12^{13} +  1 $ &$= 19^{13}$  &  $4^{13} +  1 $ &$= 11^{13}$ \\  
    $3^{9} +  1 $ &$= 11^{9}$  &  $12^{9} +  1 $ &$= 20^{9}$  &  $8^{13} +  1 $ &$= 21^{13}$  &  $2^{13} +  1 $ &$= 15^{13}$ \\  
    $4^{9} +  1 $ &$= 15^{9}$  &  $8^{9} +  1 $ &$= 19^{9}$  &  $6^{13} +  1 $ &$= 20^{13}$  &  $3^{13} +  1 $ &$= 17^{13}$ \\  
    $5^{9} +  1 $ &$= 18^{9}$  &  $5^{9} +  1 $ &$= 18^{9}$  &  $14^{13} +  1 $ &$= 9^{13}$  &  $14^{13} +  1 $ &$= 9^{13}$ \\  
    $6^{9} +  1 $ &$= 21^{9}$  &  $2^{9} +  1 $ &$= 17^{9}$  &  $4^{13} +  1 $ &$= 11^{13}$  &  $12^{13} +  1 $ &$= 19^{13}$ \\  
    $7^{9} +  1 $ &$= 6^{9}$  &  $17^{9} +  1 $ &$= 16^{9}$  &  $10^{13} +  1 $ &$= 4^{13}$  &  $19^{13} +  1 $ &$= 13^{13}$ \\  
    $8^{9} +  1 $ &$= 19^{9}$  &  $4^{9} +  1 $ &$= 15^{9}$  &  $3^{13} +  1 $ &$= 17^{13}$  &  $6^{13} +  1 $ &$= 20^{13}$ \\  
    $9^{9} +  1 $ &$= 13^{9}$  &  $10^{9} +  1 $ &$= 14^{9}$  &  $18^{13} +  1 $ &$= 16^{13}$  &  $7^{13} +  1 $ &$= 5^{13}$ \\  
    $10^{9} +  1 $ &$= 14^{9}$  &  $9^{9} +  1 $ &$= 13^{9}$  &  $7^{13} +  1 $ &$= 5^{13}$  &  $18^{13} +  1 $ &$= 16^{13}$ \\  
    $11^{9} +  1 $ &$= 10^{9}$  &  $13^{9} +  1 $ &$= 12^{9}$  &  $21^{13} +  1 $ &$= 7^{13}$  &  $16^{13} +  1 $ &$= 2^{13}$ \\ 
\hline
\end{tabular}
\end{table}

\clearpage

\section{Conclusion}
This work provides new theoretical results about the nature of the solutions of the Diophantine equation $x^{n} + y^{n} = z^{n}$ in the finite fields $\zp$ 
and their relationship with the prime $p$.  It also describes an algorithmic means of constructing all such solutions. The algorithm is not optimal, it can clearly be made more efficient, yet it works reasonably fast for prime numbers which are not too big.\\

It would interesting to extend this work to a more general setting.
\section{References}
[1] Burton, David: Elementary Number Theory (seventh edition) McGraw Hill, 2011. 
[2] Lang, Serge: Algebra (third edition) Addison Wesley, 1995. 
[3] Jacobson, Nathan: Basic Algebra I by W. H. Freeman and Company, 1974. 
[4] Herstein i.n.: Topics in Algebra (second edition) Xerox College Publishing, 1975. 
[6] Adler, Andrew and Coury, John: The Theory of Numbers by Jones and Bartlett Publishers, 1995.

\end{document}